\documentclass[a4paper,12pt]{article}

\usepackage{amsmath}
\usepackage{amsthm}
\usepackage{amssymb}  
\usepackage{latexsym} 
\usepackage[all]{xy}
\usepackage{comment}
\usepackage{rotating}
\usepackage{slashbox}
\usepackage{multirow}
\usepackage{rotating}
\usepackage[french,english]{babel}

\newcommand{\Fp}{{\mathbb{F}_p}}

\newcommand{\GL}{\operatorname{GL}}

\newcommand{\isom}{ \cong }

\newcommand{\Gal}{\operatorname{Gal}}

\newcommand{\PP}{{\mathbb P}}
\newcommand{\Q}{{\mathbb Q}}

\newcommand{\R}{{\mathbb R}}

\newcommand{\Reg}{\operatorname{Reg}}

\newcommand{\Frob}{\operatorname{Frob}}

\newcommand{\ord}{\operatorname{ord}}

\newcommand{\tor}{\operatorname{tor}}

\newcommand{\Z}{{\mathbb Z}}

\newfont{\wncyr}{wncyr10 at 12pt}
\newfont{\wncyrten}{wncyr10 at 10pt}

\newenvironment{Proof}{\par\noindent{\sc Proof:}}%
                      {\hspace*{\fill}\nobreak$\Box$\par\medskip}
                       {\hspace*{\fill}\nobreak$\Box$\par\medskip}
\newenvironment{myitemize}
{\begin{itemize}
\setlength{\itemsep}{1pt}
\setlength{\parskip}{0pt}
\setlength{\parsep}{0pt}}
{\end{itemize}}

\newtheorem{Proposition}{Proposition}[section]
\newtheorem{Theorem}[Proposition]{Theorem}
\newtheorem{Lemma}[Proposition]{Lemma}
\newtheorem{Corollary}[Proposition]{Corollary}

\theoremstyle{definition}

\newtheorem{Remark}[Proposition]{Remark}

\newcounter{nootje}
\begin{document}
\normalsize
\title{On dependence of rational points on elliptic curves}
\author{Mohammad Sadek}
\date{}
\maketitle
\let\thefootnote\relax
\begin{abstract}{\footnotesize Let $E$ be an elliptic curve defined over $\Q$. Let $\Gamma$ be a subgroup of $E(\Q)$ and $P\in E(\Q)$. In \cite{Arithmetic}, it was proved that if $E$ has no nontrivial rational torsion points, then $P\in\Gamma$ if and only if $P\in \Gamma$ mod $p$ for finitely many primes $p$. In this note, assuming the General Riemann Hypothesis, we provide an explicit upper bound on these primes when $E$ does not have complex multiplication and either $E$ is a semistable curve or $E$ has no exceptional prime. }
\end{abstract}
\begin{otherlanguage}{french}
\begin{abstract}
{\footnotesize   Soit $E$ une courbe elliptique d\'{e}finie sur $\Q$. Soit $ \Gamma$ un sous-groupe de $ E(\Q) $ et $ P \in E (\Q) $. Dans \cite{Arithmetic}, il a \'{e}t\'{e} prouv\'{e} que si $ E $ n'a pas de points de torsion rationnelles non triviales, alors $ P \in \Gamma $ si et seulement si $ P \in \Gamma $ mod $ p $ pour un nombre fini de nombres premiers $ p $. Dans cette note, supposons l'hypoth\`{e}se g\'{e}n\'{e}ral de Riemann, nous fournissons une borne-sup\'{e}rieure explicite sur ces nombres premiers quand $ E $ n'a pas de multiplication complexe et soit $ E $ est une courbe semi-stable ou $ E $ n'a aucun nombre premier exceptionnel.                 }
\end{abstract}
\end{otherlanguage}

\section{Introduction}

 For a given set of rational points on an elliptic curve $E$ defined over $\Q$, there are several methods to check if these points are linearly dependent. These methods include heights on elliptic curves, and the two descent algorithm. In \cite{Arithmetic}, the authors showed that linear dependence of rational points on certain abelian varieties over a given number field $K$ satisfies a local to global principle. Namely, a set of rational points on such an abelian variety satisfies a dependence relation over $K$ if and only if it satisfies a dependence relation when reduced modulo all but finitely many primes of $K$. In fact, they even proved a stronger version of the latter result. More precisely, a dependence relation of rational points holds on these abelian varieties if and only if these points satisfy dependence relations modulo finitely many primes. The reader interested in detecting dependence of rational points on abelian varieties via reduction maps may consult \cite{GAJDA, PJOSSEN} and the references there.

In this note, we analyse the aforementioned results. Given an elliptic curve $E/\Q$ and a basis $P_1,\ldots,P_r$ for $E(\Q)$, a point $P\in E(\Q)$ lies in a subgroup $\Gamma\subset E(\Q)$ if and only if the reduction of $P$ lies in $\Gamma$ modulo finitely many primes. The choice of these primes depends on $E$, the points $P,P_1,\ldots,P_r$, and the subgroup $\Gamma$. Assuming that $E$ has no nontrivial rational torsion points, we introduce an explicit upper bound on these primes when $E$ has no complex multiplication and either $E$ is semistable or $E$ has no exceptional prime where certain values are not attained by the $j$-invariant of $E$. The key idea in order to provide such a bound is to use an effective version of Chebotarev's theorem which assumes the Generalized Riemann hypothesis, GRH.

\section{Linear dependence of rational points}
In this section we will review the main results of \cite{Arithmetic} for an elliptic curve defined over $\Q$.
 Let $E$ be an elliptic curve defined over $\Q$. We assume moreover that $E$ has no nontrivial torsion over $\Q$. Let $P_1,\ldots,P_r$ be a basis for $E(\Q)$.

 For each $j,1\le j\le r$, a lattice $\widetilde{\Gamma}_{j}\subset E(\Q)$ is defined, see p. 334 for the precise definition of $\widetilde{\Gamma}_j$. Given $P\in E(\Q)$ there exist $n_1,\ldots,n_r\in\Z$ such that $P=n_1P_1+\ldots+n_rP_r$.
 For each $1\le j\le r$ and for each prime $l\mid n_j$, the field $L_{j,l}$ is defined as follows:
\[L_{j,l}:=\Q\left(E[l^{k_{j,l}+1}],\frac{1}{l^{k_{j,l}}}\widetilde{\Gamma}_j\right),\]
where $k_{j,l}$ is chosen such that the image of the residual representation
\[\overline{\rho}_{l^{k_{j,l}+1}}:\Gal\left(\Q(E[l^{k_{j,l}}])/\Q\right)\to  \GL_2\left(\Z/l^{k_{j,l}+1}\Z\right)\]
contains a nontrivial homothety and such that $E[l^{k_{j,l}}]$ contains at least $r$ points, see Theorem 6.3 of \cite{Arithmetic}. In particular since $E[l^{k_{j,l}}]\isom \Z/l^{k_{j,l}}\Z\times \Z/l^{k_{j,l}}\Z$, it follows that $l^{k_{j,l}}\ge\sqrt{r}$. Moreover one must have $k_{j,l}\ge \ord_l(n_j)$, see Theorem 6.4 in \cite{Arithmetic}. Therefore one may choose $\displaystyle k_{j,l}\ge \max\left\{\frac{\log r}{2 \log l},\frac{\log |n_j|}{\log l}\right\}$ such that the image of the representation $\overline{\rho}_{l^{k_{j,l}+1}}$ contains a nontrivial homothety.

One observes that \[ L_{j,l}\subset F_{j,l}:= \Q\left(E[l^{k_{j,l}+1}],\frac{1}{l^{k_{j,l}}}E(\Q)\right)\] since $\widetilde{\Gamma}_{j}\subset E(\Q)$.

The authors in \cite{Arithmetic} defined a finite set $S_{j,l}$ which contains all primes $q$ such that every $\sigma\in\Gal(L_{j,l}/\Q)$ is a Frobenius element at some $q\in S_{j,l}$. An effective version of Chebotarev's theorem was used to construct the set $S_{j,l}$.

Given a number field $L$, an effective version of Chebotarev's theorem by Lagarias and Odlyzko states that
there are effectively computable constants $b_{1,L}$ and $b_{2,L}$ such that every element $\sigma\in \Gal(L/\Q)$ is equal to the Frobenius element $\Frob_{q}\in\Gal(L/\Q)$ for an integer prime $q$ such that $q \le b_{1,L} \Delta_{L}^{b_{2,L}}$ where $\Delta_{L}$ is the discriminant of $L$.

For each $j$ such that $n_j\ne 0$, the following sets were defined in \cite{Arithmetic}
\begin{eqnarray*}S_{j,l}&:=&\{q:q\le b_{1,L_{j,l}}\Delta_{L_{j,l}}^{b_{2,L_{j,l}}}\textrm{ and $q$ is a good prime for $E$}\},\\
S_j&:=& \bigcup_{l\mid n_j} S_{j,l}.
\end{eqnarray*}
The set $S$ is defined by
\[S:=\bigcup_{1\le j\le r,n_j\ne 0}S_j.\]

A local to global property for dependence of rational points on an abelian variety of certain type defined over a number field can be found in \cite{Arithmetic} and the references there. In fact, an elliptic curve defined over $\Q$ is an abelian variety which satisfies the hypotheses of \cite[Theorem 6.4]{Arithmetic}, see \cite[Corollary 4.3]{Arithmetic}. Throughout this note, if $P\in\PP^2(\Q)$, we write $P_p$ for the reduction of $P$ in $\PP^2(\Fp)$.
\begin{Theorem}[Theorem 6.4, \cite{Arithmetic}]
\label{thm1}
Let $E$ be an elliptic curve defined over $\Q$. Let $P\in E(\Q)$ and let $\Gamma$ be a subgroup of $E(\Q)$. Let $S$ be the finite set defined above. If $P_p \in \Gamma$ mod $p$
for all $p\in S$ then $P\in\Gamma + E(\Q)_{\tor}$.
Hence if $E(\Q)_{\tor}\subset \Gamma$, then the following conditions are equivalent:
\begin{myitemize}
\item[(1)] $P \in\Gamma$.
\item[(2)] $P_p \in \left(\Gamma\textrm{ mod }p\right)$ for all $p\in S$.
\end{myitemize}
\end{Theorem}

\begin{Remark}
\label{rem1}
In Theorem \ref{thm1}, the equivalence holds because $S_{j,l}$ contains the primes $q$ for which every $\sigma\in\Gal(L_{j,l}/\Q)$ is equal to the Frobenius element $\Frob_q$ at $q$. An effective version of Chebotarev's theorem by Lagarias and Odlyzko is used to provide an upper bound for these primes $q$, see the proofs of Theorem 6.2 and Theorem 6.4 of \cite{Arithmetic}. Thus one can replace $S_{j,l}$ with any finite set containing the primes $q$ in the definition of $S$ in Theorem \ref{thm1}. In fact, we will use a different effective version of Chebotarev's theorem to introduce an alternative finite set.
\end{Remark}


In the following lemma, we collect different effective versions of Chebotarev's Density Theorem. One can use these versions to redefine the sets $S_{j,l},S_j$ and $S$, see Remark \ref{rem1}. The following can be found as Theorem 2.2 and Proposition 2.3 in \cite{Lichtenstein}.

\begin{Lemma}
\label{lem2}
Let $L/\Q$ be a finite Galois extension. We denote the absolute value of the discriminant and the degree of $L/\Q$ by $\Delta_L$ and $d_L$ respectively. Let $C$ be a conjugacy class of $\Gal(L/\Q)$. There is an integer prime $p$ such that the Frobenius at $p$ is in $C$, and such that $p$ satisfies the following bounds.
\begin{myitemize}
\item[a)] There is an absolute effectively computable constant $A$ such that $p\le 2\Delta_L^A$.
\end{myitemize}
Now we assume the GRH.
\begin{myitemize}
\item[b)] There is an absolute effectively computable constant $b$ such that $p\le b(\log \Delta_L)^2$. In fact, one may take $b=70$.
\item[c)] If $S$ is a set of prime numbers such that $L/\Q$ is unramified outside of $S$, for the conjugacy class $C$ in $\Gal(L/\Q)$, there exists a prime number $p\not\in S$ such that the Frobenius at $p$ is in $C$, and such that \[p\le 280d_L^2\left(\log d_L+\sum_{q\in S}\log q\right)^2. \]
\end{myitemize}
\end{Lemma}

For each $1\le j\le r$ and each prime $l\mid n_j$, we recall that $\displaystyle L_{j,l}\subset F_{j,l}:= \Q\left(E[l^{k_{j,l}+1}],\frac{1}{l^{k_{j,l}}}E(\Q)\right)$.
Now we define sets $S'_{j,l},S'_j$ and $S'$ under the assumption of the GRH. Assuming the GRH, one can use Lemma \ref{lem2} (c) in order to define a set $S'_{j,l}$ which contains all the primes $q$ such that every $\sigma\in\Gal(F_{j,l}/\Q)$ is the Frobenius element at some $q\in S'_{j,l}$. We set
\begin{eqnarray*}
S'_{j,l}&:=&\{q:q\le 280d_{F_{j,l}}^2\left(\log d_{F_{j,l}}+\sum_{q\in B}\log q\right)^2\textrm{ and $q$ is a good prime for $E$}\},\\
S'_j&:=& \bigcup_{l\mid n_j} S'_{j,l},\\
S'&:=&\bigcup_{1\le j\le r,n_j\ne 0}S'_j,
\end{eqnarray*}
where $d_{F_{j,l}}$ is the degree of the extension $F_{j,l}/\Q$, and $B$ is the set of primes outside which the field $F_{j,l}$ is unramified. In fact, the field $F_{j,l}$ is unramified outside the set of bad primes of $E$ and the prime $l$, see \cite[Theorem 1]{Serre-Tate}.

Theorem \ref{thm1} and Remark \ref{rem1} yield the following consequence.
\begin{Corollary}
\label{cor1}
 Assume that $E$ is an elliptic curve defined over $\Q$ such that $E(\Q)_{\tor}=\{O_E\}$. Let $P\in E(\Q)$ and $\Gamma$ a subgroup of $E(\Q)$. Let $S'$ be defined as above. The following conditions are equivalent:
\begin{myitemize}
\item[(1)] $P \in\Gamma$.
\item[(2)] $P_p \in \left(\Gamma\textrm{ mod }p\right)$ for all $p\in S'$.
\end{myitemize}
\end{Corollary}

\section{Bounds}
In this section we find explicit bounds for the coefficients of a linear dependence relation in $E(\Q)$.

\subsection{A bound on the coefficients of a linear dependence relation}

Let $E$ be an elliptic curve defined over $\Q$ with rank $r$ such that $E(\Q)_{\tor}=\{O_E\}$. Let $P_1,\ldots,P_r$ be a basis for $E(\Q)$. Recall that the height pairing on $E$ is
\begin{eqnarray*}
\langle\textrm{\hskip4pt},\textrm{\hskip4pt}\rangle&:& E\times E\to \R\\
\langle P,Q\rangle&=&\hat{h}(P+Q)-\hat{h}(P)-\hat{h}(Q)
\end{eqnarray*}
 where $\hat{h}:E\to \R$ is the canonical height on $E$, see (\S 9, VIII, \cite{sil1}). The regulator matrix $R_E$ of $E$ is given by $\displaystyle \left(\langle P_i,P_j\rangle\right)_{1\le i,j\le r}$. The eigen values of $\Reg_E$ are $\lambda_1\le\lambda_2\le\ldots\le \lambda_r$.

 \begin{Lemma}
 \label{lem:Reg}
 Let $P\in E(\Q)$ be such that
\[P=\sum_{i=1}^rn_iP_i,\;n_i\in\Z.\] For each $j$, $1\le j\le r$, one has \[|n_j|\le \left|\frac{\langle P,P\rangle}{\lambda_1}\right|^{1/2}\] where $\displaystyle N^T=\left(
                         \begin{array}{cccc}
                           n_1 & n_2 & \ldots & n_r \\
                         \end{array}
                       \right)
$.
\end{Lemma}
\begin{Proof}
Using the fact that the height pairing $\langle\textrm{\hskip4pt},\textrm{\hskip4pt}\rangle$ is bilinear and positive definite, one obtains the following equality \[\langle P,P\rangle=\sum_{i,j} n_i n_j \langle P_i,P_j\rangle=N^T\Reg_E N\ge 0.\] In addition, one has \[\min_i\lambda_i\le \frac{N^T \Reg_E N}{N^TN}\le \max_i\lambda_i.\]
One observes that $\displaystyle N^TN=\sum_i n_i^2$, therefore for any $j$, $1\le j\le r$, one has
 \[\langle P,P\rangle=N^T\Reg_E N\ge \lambda_1 N^TN\ge \lambda_1n_j^2.\]
\end{Proof}

\subsection{Bounding the degree of a Kummer extension}

In the following lemma we estimate the degree of the Kummer extension $\displaystyle F_{j,l}=\Q\left(E[l^{k_{j,l}+1}],\frac{1}{l^{k_{j,l}}}E(\Q)\right)$.

\begin{Lemma}
\label{lem3}
Let $E/\Q$ be an elliptic curve of rank $r$. Let $l$ be an integer prime and $m$ a positive integer. Let $L$ denote the field $\displaystyle\Q\left(E[l^{m}],\frac{1}{l^{m-1}}E(\Q)\right)$. Then \[[L:\Q]\le (l^2 - 1)(l^2 - l)l^{2mr+4(m-1)}.\]
\end{Lemma}
\begin{Proof}
The Galois group of the field extension $\displaystyle \Q\left(E[l^{m}],\frac{1}{l^{m-1}}E(\Q)\right)/\Q(E[l^{m}])$ can be viewed as a subgroup of the product $\displaystyle \left(E[l^{m}]\right)^r$. Therefore, the degree of the extension is at most $l^{2mr}$. Now it is known that $\displaystyle\Gal\left(\Q(E[l^{m}])/\Q\right)\hookrightarrow\GL_2(\Z/l^{m}\Z)$, where the latter group is of order $(l^2 - 1)(l^2 - l)l^{4(m-1)}$, hence follows the upper bound for $[L:\Q]$.
\end{Proof}

\section{Main results}

We assume that $E$ is an elliptic curve defined over $\Q$ with $E(\Q)_{\tor}=\{O_E\}$ and the rank of $E(\Q)$ is $r>0$. Therefore, all rational points in $E(\Q)$ are of infinite order. Let $P_1,\ldots,P_r$ be a basis for $E(\Q)$.

  Let $\Gamma=d_1\Z P_1+d_2\Z P_2+\ldots+d_s\Z P_s+\ldots+d_r\Z P_r$ where $d_i\in\Z\setminus\{0\}$ if $1\le i\le s$ and $d_i=0$ otherwise. Assume that $P=n_1P_1+\ldots+n_r P_r$ for some $n_i\in \Z$. In view of Corollary \ref{cor1}, $P\in\Gamma$ if and only if $P\in \Gamma$ mod $p$ for every prime $p$ lying in the finite set $S'$.

  \begin{Lemma}
  \label{lem4}
  Assume the GRH. For any prime $q\in S'_{j,l}$, one has \[q\le 280M_{j,l}^2\left[\log (M_{j,l}\,Q_{j,l} )\right]^2\] where $M_{j,l}=(l^2 - 1)(l^2 - l)l^{2(k_{j,l}+1)r+4k_{j,l}}$, and $Q_{j,l}$ is the product of the prime divisors of the discriminant of the field extension $F_{j,l}/\Q$.
  \end{Lemma}
  \begin{Proof}
  Recall that \[S'_{j,l}:=\{q:q\le 280d_{F_{j,l}}^2\left(\log d_{F_{j,l}}+\sum_{p\in B_{j,l}}\log p\right)^2\textrm{ and $q$ is a good prime for $E$}\}\] where $B_{j,l}$ is the set of primes outside which $F_{j,l}$ is unramified. Those primes are exactly the prime divisors of the discriminant of the field extension $F_{j,l}/\Q$. In particular, the set $B_{j,l}$ contains the bad primes of $E$ together with $l$. We set $d_{F_{j,l}}$ to be the degree of the extension $F_{j,l}/\Q$.  Therefore for every $q\in S'_{j,l}$, one has $\displaystyle q\le 280d_{F_{j,l}}^2\left(\log d_{F_{j,l}}+\sum_{l'} \log l'\right)^2$ where $l'$ is a prime in $B_{j,l}$. Now the statement of the lemma follows once one observes that $M_{j,l}$ is the upper bound of $d_{F_{j,l}}$ obtained in Lemma \ref{lem3}.
  \end{Proof}

  We recall that a prime integer $l$ is said to be an {\em exceptional prime} for an elliptic curve $E$ defined over $\Q$ if the mod $l$ Galois representation $\rho_{E,l}:\Gal\left(\Q(E[l])/\Q\right)\to\GL_2(\Z/l\Z)$ is not surjective. If $E$ has complex multiplication then every prime is exceptional except possibly for the prime $2$. If $E$ has no complex multiplication then it was proved by Serre that the number of exceptional primes is finite, see \cite{Serre}. In fact Serre conjectured that any exceptional prime for $E$ is less than or equal to $37$. Mazur proved that if $E$ is semistable with no complex multiplication, then no prime $\ge 11$ can be exceptional for $E$, see \cite{Mazur}. In \cite{Duke}, in terms of heights, almost all elliptic curves are proved to have no exceptional primes.

  We will use the following two lemmas to produce an explicit bound on the primes in the finite set $S'$ defined in \S 2 if $E$ defined over $\Q$ has no complex multiplication, and either $E$ has no exceptional primes or $E$ is semistable.

  \begin{Lemma}
  \label{lem:GaloisSurjectivity1}
  Let $E$ be an elliptic curve defined over $\Q$ with no complex multiplication. Let $j$ be the $j$-invariant of $E$. Let $\rho_{l^n}:\Gal\left(\Q(E[l^n])/\Q\right)\to \GL_2(\Z/l^n\Z)$ be the Galois representation associated to the ${l^n}$-torsion points of $E$. The following statements hold.
  \begin{myitemize}
  \item[i.] $\rho_2$ is not surjective if and only if $j=256(t+1)^3/t$ or $j=t^2+1728$ for some $t\in\Q$.
  \item[ii.] $\rho_3$ is not surjective if and only if $j = 27(t + 1)(t + 9)^3/t^3$ or $j = t^3$ for some $t\in\Q$.
\item[iii.]  $\rho_5$ is not surjective if and only if \[ j=\frac{5^3(t + 1)(2t + 1)^3(2t^2 -3t + 3)^3}{(t^2 + t - 1)^5},\; j=\frac{5^2(t^2 + 10t + 5)^3}{t^5},\;\textrm{ or }j= t^3(t^2 + 5t + 40)\]
\item[iv.] $\rho_7$ is not surjective if and only if
\begin{eqnarray*}
j &=\frac{t(t + 1)^3(t^2 - 5t + 1)^3(t^2 - 5t + 8)^3(t^4 - 5t^3 + 8t^2 - 7t + 7)^3}{(t^3 -4t^2 + 3t + 1)^7},\\
j &=\frac{64t^3(t^2 + 7)^3(t^2 - 7t + 14)^3(5t^2 - 14t - 7)^3}{(t^3 - 7t^2 + 7t + 7)^7},\textrm{ or }\\
j &=\frac{(t^2 + 245t + 2401)^3(t^2 + 13t + 49)}{t^7}
\end{eqnarray*} for some $t\in\Q$.
  \end{myitemize}
  \end{Lemma}
  \begin{Proof}
   This is Proposition 6.1 in \cite{Zywina}.
  \end{Proof}
  \begin{Lemma}
  \label{lem:GaloisSurjectivity2}
  Let $E$ be an elliptic curve defined over $\Q$ with no complex multiplication. Let $j$ be the $j$-invariant of $E$. Let $\rho_{l^{\infty}}:\Gal\left(\overline{\Q}/\Q\right)\to \GL_2(\Z_l)$ be the Galois representation describing the Galois action on the Tate module of $E$. The following statements hold.
  \begin{myitemize}
  \item[i.] The representation $\rho_{2^{\infty}}$ is not surjective if and only if $\rho_2$ is not surjective, or $j$ is of the form
  \[-4t^3(t+8),\;-t^2+1728,\;2t^2+1728,\;\textrm{or }-2t^2+1728\] for some $t\in\Q$.
  \item[ii.] The representation $\rho_{3^{\infty}}$ is not surjective if and only if $\rho_3$ is not surjective, or
  \[j= -\frac{3^7(t^2 - 1)^3(t^6 + 3t^5 + 6t^4 + t^3 - 3t^2 + 12t + 16)^3(2t^3 + 3t^2 - 3t - 5)}{(t^3 - 3t - 1)^9}\]
  for some $t\in \Q$.
  \item[iii.] If $l\ge 5$ then $\rho_{l^{\infty}}$ is not surjective if and only if $\rho_l$ is not surjective.
  \end{myitemize}
  \end{Lemma}
  \begin{Proof}
  This is Lemma 6.6 in \cite{Zywina}.
  \end{Proof}
  One remarks that if $\rho_{l^{\infty}}$ is surjective then $\rho_{l^n}$ is surjective for any positive integer $n$.
  \begin{Theorem}
  \label{thm:exceptional}
   We assume the GRH. Let $E/\Q$ be an elliptic curve with no complex multiplication and no exceptional primes. Assume that $E(\Q)\isom\Z^{r},r>0$. Let $P_1,\ldots,P_r$ be a basis for $E(\Q)$, and $\lambda$ the minimum eigen value of the regulator matrix $\displaystyle(\langle P_i,P_j\rangle)_{1\le i,j\le r}$.
   We assume that the $j$-invariant $j$ of $E$ is not written as
   \begin{eqnarray*} \frac{256(t+1)^3}{t}, \;t^2+1728, \;-4t^3(t+8),\;-t^2+1728,\;2t^2+1728,\;\textrm{or }-2t^2+1728
   \end{eqnarray*}
   for any $t\in\Q$. Let $P\in E(\Q)$ and $\Gamma$ a subgroup of $E(\Q)$. The following conditions are equivalent.
 \begin{myitemize}
 \item[1)] $P\in \Gamma$ where $P=n_1P_1+\ldots+n_rP_r$, $n_i\in\Z$.
 \item[2)] $P_p\in\Gamma$ mod $p$, for every prime $p\le280M^2\left[\log (M\, Q )\right]^2$
 where $M=(C-1)\left(C-\sqrt{C}\right)C^{K(r+2)+r}$, $\displaystyle C=\left|\frac{\langle P,P\rangle}{\lambda}\right|^{1/2}$, $\displaystyle K=\max\left\{2,\frac{\log r}{2\log 2},\frac{\log C}{\log 2}\right\}$, and $\displaystyle Q=\max_{\substack{1\le j\le r\\l\mid n_j}}Q_{j,l}$ where $Q_{j,l}$ is the product of the prime divisors of the discriminant of $F_{j,l}/\Q$.
 \end{myitemize}
  \end{Theorem}
  \begin{Proof}
    Theorem \ref{thm1} implies that $P\in\Gamma$ if and only if $P_p\in\Gamma$ mod $p$ for every $p\in S'$ where
 \begin{eqnarray*}
 S':=\bigcup_{1\le j\le r,\,n_j\ne 0 }\left( \bigcup_{l| n_j} S'_{j,l}\right).
 \end{eqnarray*}
Lemma \ref{lem4} implies that if $q\in S'_{j,l}$, then $q\le 280M_{j,l}^2\left[\log (M_{j,l}\, Q_{j,l} )\right]^2$ where $M_{j,l}=(l^2 - 1)(l^2 - l)l^{2(k_{j,l}+1)r+4k_{j,l}}$, and we choose $\displaystyle k_{j,l}\ge \max\left\{\frac{\log r}{2 \log l},\frac{\log |n_j|}{\log l}\right\}$ such that the image of the residual representation
\[\overline{\rho}_{l^{k_{j,l}+1}}:\Gal\left(\Q(E[l^{k_{j,l}}])/\Q\right)\to  \GL_2\left(\Z/l^{k_{j,l}}\Z\right)\]
contains a nontrivial homothety, see \S 2. Since $E$ has no exceptional primes it follows that $\overline{\rho}_{l^2}$ is surjective for every prime $l$, hence the image of $\overline{\rho}_{l^2}$ contains a nontrivial homothety for every prime $l\ne 2$. It follows that one can set $\displaystyle k_{j,l}=\max\left\{1, \frac{\log r}{2 \log l},\frac{\log |n_j|}{\log l}\right\}$ when $l\ne 2$. In view of Lemma \ref{lem:GaloisSurjectivity1} and Lemma \ref{lem:GaloisSurjectivity2} our assumption on the values taken by the $j$-invariant of $E$ forces $\rho_{2^{\infty}}$ to be surjective, and hence the representation $\rho_{2^n}$ is surjective for any positive integer $n$. In particular the residual representation $\overline{\rho}_8$ is surjective. In other words the image of $\overline{\rho}_8$ contains a nontrivial homothety. Therefore one may set $\displaystyle k_{j,2}=\max\left\{2,\frac{\log r}{2 \log 2},\frac{\log |n_j|}{\log 2}\right\}$.

 For every $j$, $1\le j\le r$, one has $|n_j|\le C$, see Lemma \ref{lem:Reg}. If $l$ is a prime dividing $n_j$, then $2\le l\le \sqrt{C}$. Therefore $\displaystyle k_{j,l}\le K=\max\left\{2,\frac{\log r}{2\log 2},\frac{\log C}{\log 2}\right\}$ for every $j,l$. Thus one obtained an upper bound $M=(C-1)\left(C-\sqrt{C}\right)C^{(K+1)r+2K}$ for $M_{j,l}$ in Lemma \ref{lem4} for any $l,j$. Therefore if $p\in S'$, then $p\le 280M^2\left[\log (M\, Q )\right]^2$.
  \end{Proof}

 \begin{Theorem}
 We assume the GRH. Let $E/\Q$ be a semistable elliptic curve with no complex multiplication such that $E(\Q)\isom\Z^{r},r>0$. Let $P_1,\ldots,P_r$ be a basis for $E(\Q)$, and $\lambda$ the minimum eigen value of the regulator matrix $\displaystyle(\langle P_i,P_j\rangle)_{1\le i,j\le r}$. We assume that the $j$-invariant $j$ of $E$ is not written as
{\footnotesize \begin{eqnarray*}
\frac{256(t+1)^3}{t}, \;t^2+1728, \;-4t^3(t+8),\;-t^2+1728,\;2t^2+1728,\;-2t^2+1728,\\
  \frac{27(t + 1)(t + 9)^3}{t^3},\; t^3,\;-\frac{3^7(t^2 - 1)^3(t^6 + 3t^5 + 6t^4 + t^3 - 3t^2 + 12t + 16)^3(2t^3 + 3t^2 - 3t - 5)}{(t^3 - 3t - 1)^9},\\
  \frac{5^3(t + 1)(2t + 1)^3(2t^2 -3t + 3)^3}{(t^2 + t - 1)^5}, \frac{5^2(t^2 + 10t + 5)^3}{t^5},\; t^3(t^2 + 5t + 40),\\
  \frac{t(t + 1)^3(t^2 - 5t + 1)^3(t^2 - 5t + 8)^3(t^4 - 5t^3 + 8t^2 - 7t + 7)^3}{(t^3 -4t^2 + 3t + 1)^7},\\
\frac{64t^3(t^2 + 7)^3(t^2 - 7t + 14)^3(5t^2 - 14t - 7)^3}{(t^3 - 7t^2 + 7t + 7)^7}, \textrm{ or }
\frac{(t^2 + 245t + 2401)^3(t^2 + 13t + 49)}{t^7}
  \end{eqnarray*}}
  for any $t\in\Q$.
 Let $P\in E(\Q)$ and $\Gamma$ a subgroup of $E(\Q)$. The following conditions are equivalent.
 \begin{myitemize}
 \item[1)] $P\in \Gamma$ where $P=n_1P_1+\ldots+n_rP_r$, $n_i\in\Z$.
 \item[2)] $P_p\in\Gamma$ mod $p$, for every prime $p\le280M^2\left[\log (M\, Q )\right]^2$
 where  $M=(C-1)\left(C-\sqrt{C}\right)C^{K(r+2)+r}$, $\displaystyle C=\left|\frac{\langle P,P\rangle}{\lambda}\right|^{1/2}$, $\displaystyle K=\max\left\{2,\frac{\log r}{2\log 2},\frac{\log C}{\log 2}\right\}$, and $\displaystyle Q=\max_{\substack{1\le j\le r\\l\mid n_j}}Q_{j,l}$ where $Q_{j,l}$ is the product of the prime divisors of the discriminant of $F_{j,l}/\Q$.
 \end{myitemize}
 \end{Theorem}
 \begin{Proof}
 The proof is similar to that of Theorem \ref{thm:exceptional}. Since $E$ is semistable with no complex multiplication, it follows that the residual representation
\[\overline{\rho}_{l^2}:\Gal\left(\Q(E[l])/\Q\right)\to  \GL_2\left(\Z/l\Z\right) \]
is surjective for any $l\ge 11$, see \cite{Serre}. That the $j$-invariant is not one of the ones above implies that $\overline{\rho}_{l^2}$ is surjective for any prime $l<11$, see Lemma \ref{lem:GaloisSurjectivity1} and Lemma \ref{lem:GaloisSurjectivity2}. It follows that the image of $\overline{\rho}_{l^2}$ contains a nontrivial homothety for any $l\ne 2$. Again our assumption on the $j$-invariant implies that $\rho_{2^{\infty}}$ is surjective which yields that
 the residual representation $\overline{\rho}_{8}$ contains a nontrivial homothety. The latter argument together with the fact that $2\le l\le \sqrt{C}$ yield the upper bound $K$ for $k_{j,l}$, hence the upper bound $M$ for $M_{j,l}$ for every $j,l$, see Lemma \ref{lem4}.
 \end{Proof}

 \hskip-18pt\emph{\bf{Acknowledgements.}}
I would like to thank the referee for many comments, corrections,
suggestions, and a great deal of patience that helped the author improve the manuscript significantly.
\bibliographystyle{plain}
\footnotesize
\bibliography{ECDLP}
Department of Mathematics and Actuarial Science\\ American University in Cairo\\ mmsadek@aucegypt.edu
\end{document}